# Similar Ruled Surfaces with Variable Transformations in Minkowski 3-space $E_1^3$

**Mehmet Önder**
*Celal Bayar University, Faculty of Science and Arts, Department of Mathematics, Muradiye Campus, 45047 Muradiye, Manisa, Turkey. E-mail: mehmet.onder@cbu.edu.tr*

**Abstract**
In this study, we consider the notion of similar ruled surface for timelike and spacelike ruled surfaces in Minkowski 3-space $E_1^3$. We obtain some properties of these special surfaces in $E_1^3$ and we show that developable ruled surfaces in $E_1^3$ form a family of similar ruled surfaces if and only if the striction curves of the surfaces are similar curves with variable transformation. Moreover, we obtain that cylindrical surfaces and conoids form two families of similar ruled surfaces in $E_1^3$.

## 1. Introduction

In the curve theory, special curve pairs for which at the corresponding points of the curves one of the Frenet vectors of a curve coincides with one of the Frenet vectors of other curve, are very interesting and an important problem of the differential geometry. Bertrand curves, Mannheim curves and involute-evolute curves are the well-known of such curves and studied extensively [5,12,16]. Recently, a new definition of the special curves was given by El-Sabbagh and Ali [2]. They have called these new curves as similar curves with variable transformation and defined as follows: Let $\psi_\alpha(s_\alpha)$ and $\psi_\beta(s_\beta)$ be two regular curves in $E^3$ parameterized by arc lengths $s_\alpha$ and $s_\beta$ with curvatures $\kappa_\alpha$, $\kappa_\beta$ and torsions $\tau_\alpha$, $\tau_\beta$ and Frenet frames $\{\vec{T}_\alpha, \vec{N}_\alpha, \vec{B}_\alpha\}$ and $\{\vec{T}_\beta, \vec{N}_\beta, \vec{B}_\beta\}$. $\psi_\alpha(s_\alpha)$ and $\psi_\beta(s_\beta)$ are called similar curves with variable transformation $\lambda_\beta^\alpha$ if there exists a variable transformation

$$s_\alpha = \int \lambda_\beta^\alpha(s_\beta) ds_\beta,$$

of the arc lengths such that the tangent vectors are the same for the two curves i.e., $\vec{T}_\alpha = \vec{T}_\beta$ for all corresponding values of parameters under the transformation $\lambda_\beta^\alpha$. They have called all curves satisfying this condition as a family of similar curves. Moreover, they have obtained some properties of the family of similar curves.

Furthermore, the surface pairs especially ruled surface pairs (called offset surfaces) have an important positions and applications in the study of design problems in spatial mechanisms and physics, kinematics and computer aided design (CAD) [9,10]. So, these surfaces are one of the most important topics of the surface theory. In fact, ruled surface offsets are the generalizations of the notion of Bertrand curves, Mannheim curves and similar curves to the line geometry and these surface pairs are called Bertrand offsets, Mannheim offsets and similar ruled surfaces, respectively [4,7,8,11].

In this work, we introduce timelike and spacelike similar ruled surfaces in Minkowski 3-space $E_1^3$. We give some theorems characterizing these special surfaces and we show that developable ruled surfaces in $E_1^3$ form a family of similar ruled surfaces if and only if the striction curves of the surfaces are similar curves with variable transformation.

## 2. Preliminaries

Let $E_1^3$ be a Minkowski 3-space with natural Lorentz Metric

$$\langle,\rangle = -dx_1^2 + dx_2^2 + dx_3^2,$$

where $(x_1, x_2, x_3)$ is a rectangular coordinate system of $E_1^3$. According to this metric, in $E_1^3$ an arbitrary vector $\vec{v} = (v_1, v_2, v_3)$ can have one of three Lorentzian causal characters; it can be spacelike if $\langle \vec{v}, \vec{v} \rangle > 0$ or $\vec{v} = 0$, timelike if $\langle \vec{v}, \vec{v} \rangle < 0$ and null (lightlike) if $\langle \vec{v}, \vec{v} \rangle = 0$ and $\vec{v} \neq 0$ [6]. Similarly, an arbitrary curve $\vec{\alpha} = \vec{\alpha}(s)$ can locally be spacelike, timelike or null (lightlike), if all of its velocity vectors $\vec{\alpha}'(s)$ are spacelike, timelike or null (lightlike), respectively. For the vectors $\vec{x} = (x_1, x_2, x_3)$ and $\vec{y} = (y_1, y_2, y_3)$ in $E_1^3$, the vector product of $\vec{x}$ and $\vec{y}$ is defined by

$$\vec{x} \wedge \vec{y} = (x_2 y_3 - x_3 y_2, x_1 y_3 - x_3 y_1, x_2 y_1 - x_1 y_2).$$

The Lorentzian sphere and hyperbolic sphere of radius $r$ and center 0 in $E_1^3$ are given by

$$S_1^2 = \{\vec{x} = (x_1, x_2, x_3) \in E_1^3 : \langle \vec{x}, \vec{x} \rangle = r^2\},$$

and

$$H_0^2 = \{\vec{x} = (x_1, x_2, x_3) \in E_1^3 : \langle \vec{x}, \vec{x} \rangle = -r^2\},$$

respectively [14,15].

Analogue to the curves, a surface can be timelike or spacelike in $E_1^3$. A surface in the Minkowski 3-space $E_1^3$ is called a timelike surface if the induced metric on the surface is a Lorentz metric and is called a spacelike surface if the induced metric on the surface is a positive definite Riemannian metric, i.e., the normal vector on spacelike (timelike) surface is a timelike (spacelike) vector [1].

## 3. Timelike and Spacelike Ruled Surface in Minkowski 3-space

Let $I$ be an open interval in the real line $IR$. Let $\vec{k} = \vec{k}(u)$ be a curve in $E_1^3$ defined on $I$ and $\vec{q} = \vec{q}(u)$ be a unit direction vector of an oriented line in $E_1^3$. Then we have the following parametrization for a ruled surface $N$,

$$\vec{\varphi}(u,v) = \vec{k}(u) + v \vec{q}(u). \qquad (1)$$

The parametric $u$-curve of this surface is a straight line of the surface which is called ruling. For $v = 0$, the parametric $v$-curve of this surface is $\vec{k} = \vec{k}(u)$ which is called base curve or generating curve of the surface. In particular, if the direction of $\vec{q}$ is constant, the ruled surface is said to be cylindrical, and non-cylindrical otherwise.

The distribution parameter (or drall) of the ruled surface in (1) is given as

$$\delta_\varphi = \frac{\left| d\vec{k}, \vec{q}, d\vec{q} \right|}{\langle d\vec{q}, d\vec{q} \rangle} \qquad (2)$$

([3]). Then the normal vectors are collinear at all points of same ruling and at nonsingular points of the surface $N$, the tangent planes are identical. We then say that tangent plane contacts the surface along a ruling. Such a ruling is called a torsal ruling. If $\left| d\vec{k}, \vec{q}, d\vec{q} \right| \neq 0$, then the tangent planes of the surface $N$ are distinct at all points of same ruling which is called nontorsal [14,15].

For the unit normal vector $\vec{m}$ of the ruled surface $N$ we have $\vec{m} = \dfrac{\vec{\varphi}_u \times \vec{\varphi}_v}{\|\vec{\varphi}_u \times \vec{\varphi}_v\|}$.

So, at the points of a nontorsal ruling $u = u_1$ we have

$$\vec{a} = \lim_{v \to \infty} \vec{m}(u_1, v) = \frac{d\vec{q} \times \vec{q}}{\|d\vec{q}\|}.$$

The point at which the unit normal of $N$ is perpendicular to $\vec{a}$ is called the striction point (or central point) $C$ and the set of striction points of all rulings is called striction curve of the surface. The parametrization of the striction curve $\vec{c} = \vec{c}(u)$ on a ruled surface is given by

$$\vec{c}(u) = \vec{k}(u) - \frac{\langle d\vec{q}, d\vec{k} \rangle}{\langle d\vec{q}, d\vec{q} \rangle} \vec{q}, \tag{3}$$

[13,14,15]. So that, the base curve of the ruled surface is its striction curve if and only if $\langle d\vec{q}, d\vec{k} \rangle = 0$.

The vector $\vec{h}$ defined by $\vec{h} = \pm \vec{a} \times \vec{q}$ is called central normal which is the surface normal along the striction curve. Then the orthonormal system $\{C; \vec{q}, \vec{h}, \vec{a}\}$ is called Frenet frame of the ruled surfaces $N$ where $C$ is the central point of ruling of ruled surface $N$ and $\vec{q}$, $\vec{h} = \pm \vec{a} \times \vec{q}$, $\vec{a}$ are unit vectors of ruling, central normal and central tangent, respectively.

Let now consider the ruled surface $N$. According to the Lorentzian casual characters of ruling and central normal, we can give the following classifications of the ruled surface $N$ as follows;

**i)** If the central normal vector $\vec{h}$ is spacelike and $\vec{q}$ is timelike, then the ruled surface $N$ is said to be of type $N_-$.

**ii)** If the central normal vector $\vec{h}$ and the ruling $\vec{q}$ are both spacelike, then the ruled surface $N$ is said to be of type $N_+$.

**iii)** If the central normal vector $\vec{h}$ is timelike, then the ruled surface $N$ is said to be of type $N_\times$ [14,15].

The ruled surfaces of type $N_+$ and $N_-$ are clearly timelike and the ruled surface of type $N_\times$ is spacelike. By using these classifications and taking the striction curve as the base curve the parametrization of the ruled surface $N$ can be given as follows,

$$\varphi(s, v) = \vec{c}(s) + v \vec{q}(s), \tag{4}$$

where $\langle \vec{q}, \vec{q} \rangle = \varepsilon (= \pm 1)$, $\langle \vec{h}, \vec{h} \rangle = \pm 1$ and $s$ is the arc length of the striction curve.

For the derivatives of the vectors of Frenet frame $\{C; \vec{q}, \vec{h}, \vec{a}\}$ of ruled surface $N$ with respect to the arc length $s$ of striction curve we have the followings

**i)** If the ruled surface $N$ is a timelike ruled surface then we have

$$\begin{bmatrix} d\vec{q}/ds \\ d\vec{h}/ds \\ d\vec{a}/ds \end{bmatrix} = \begin{bmatrix} 0 & k_1 & 0 \\ -\varepsilon k_1 & 0 & k_2 \\ 0 & \varepsilon k_2 & 0 \end{bmatrix} \begin{bmatrix} \vec{q} \\ \vec{h} \\ \vec{a} \end{bmatrix} \tag{5}$$

and

$$\vec{q} \times \vec{h} = \varepsilon \vec{a}, \quad \vec{h} \times \vec{a} = -\varepsilon \vec{q}, \quad \vec{a} \times \vec{q} = -\vec{h}, \tag{6}$$

[See 14].

**ii)** If the ruled surface $N$ is spacelike ruled surface then we have

$$\begin{bmatrix} d\vec{q}/ds \\ d\vec{h}/ds \\ d\vec{a}/ds \end{bmatrix} = \begin{bmatrix} 0 & k_1 & 0 \\ k_1 & 0 & k_2 \\ 0 & k_2 & 0 \end{bmatrix} \begin{bmatrix} \vec{q} \\ \vec{h} \\ \vec{a} \end{bmatrix}, \tag{7}$$

and

$$\vec{q} \times \vec{h} = -\vec{a}, \quad \vec{h} \times \vec{a} = -\vec{q}, \quad \vec{a} \times \vec{q} = \vec{h}, \tag{8}$$

[See 15].

In the equations (5) and (7), $k_1 = \dfrac{ds_1}{ds}$, $k_2 = \dfrac{ds_3}{ds}$ and $s_1$, $s_3$ are the arc lengths of the spherical curves circumscribed by the bound vectors $\vec{q}$ and $\vec{a}$, respectively. Moreover, timelike and spacelike ruled surfaces satisfying $k_1 \neq 0$, $k_2 = 0$ are called timelike and spacelike conoids in $E_1^3$, respectively [14,15].

Now, we can represent and prove the following theorems which are necessary for the following section.

**Theorem 3.1.** *Let the striction curve $\vec{c} = \vec{c}(s)$ of ruled surface $N$ be unit speed curve with same Lorentzian casual character with the ruling and let $\vec{c}(s)$ be the base curve of the surface. Then $N$ is developable if and only if the unit tangent of the striction curve is the same with the ruling along the curve.*

**Proof:** Let $N$ be a timelike ruled surface an let $s$ be arc length parameter of the striction curve. Then the unit tangent of the striction curve is given by

$$\vec{T}(s) = \frac{d\vec{c}}{ds} = (\cosh\theta)\vec{q}(s) + (\sinh\theta)\vec{a}(s),$$

where $\theta = \theta(s)$ is the angle between unit vectors $\vec{T}(s)$ and $\vec{q}(s)$ [14]. Since the striction curve is taken as base curve, then from (2) and (5) the distribution parameter of the surface $N$ is obtained as

$$d = -\frac{\sinh\theta}{k_1}.$$

Thus we have that timelike ruled surface $N$ is developable if and only if $\vec{T}(s) = \vec{q}(s)$ satisfies.

If $N$ is a spacelike ruled surface then the unit tangent of the striction curve is given by

$$\vec{T}(s) = \frac{d\vec{c}}{ds} = (\cos\theta)\vec{q}(s) + (\sin\theta)\vec{a}(s),$$

(See [15]). Then from (2) and (5) the distribution parameter of the surface $N$ is obtained as

$$d = \frac{\sin\theta}{k_1}.$$

Thus we have that spacelike ruled surface $N$ is developable if and only if $\vec{T}(s) = \vec{q}(s)$ satisfies.

**Theorem 3.2.** Let the striction curve $\vec{c} = \vec{c}(s)$ of ruled surface $N$ be unit speed i.e., $s$ is arc length parameter of $\vec{c}(s)$. Suppose that $\vec{c} = \vec{c}(\varphi)$ is another parametrization of the striction curve by the parameter $\varphi(s) = \int k_1(s) ds$. Then the ruling $\vec{q}$ satisfies a vector differential equation of third order given by

$$\begin{cases} \dfrac{d}{d\varphi}\left(\dfrac{1}{f(\varphi)} \dfrac{d^2 \vec{q}}{d\varphi^2}\right) + \varepsilon\left(\dfrac{1-f^2(\varphi)}{f(\varphi)}\right) \dfrac{d\vec{q}}{d\varphi} - \varepsilon\left(\dfrac{1}{f^2(\varphi)} \dfrac{df(\varphi)}{d\varphi}\right) \vec{q} = 0; \text{ if } N \text{ is timelike}, \\ \dfrac{d}{d\varphi}\left(\dfrac{1}{f(\varphi)} \dfrac{d^2 \vec{q}}{d\varphi^2}\right) - \left(\dfrac{1+f^2(\varphi)}{f(\varphi)}\right) \dfrac{d\vec{q}}{d\varphi} + \left(\dfrac{1}{f^2(\varphi)} \dfrac{df(\varphi)}{d\varphi}\right) \vec{q} = 0; \text{ if } N \text{ is spacelike}, \end{cases} \quad (9)$$

where $f(\varphi) = \dfrac{k_2(\varphi)}{k_1(\varphi)}$.

**Proof:** Let $N$ be a timelike ruled surface. If we write derivatives given in (5) according to $\varphi$, we have

$$\dfrac{d\vec{q}}{d\varphi} = \vec{h},$$

$$\dfrac{d\vec{h}}{d\varphi} = -\varepsilon \vec{q} + f(\varphi) \vec{a},$$

$$\dfrac{d\vec{a}}{d\varphi} = \varepsilon f(\varphi) \vec{h},$$

respectively, where $f(\varphi) = \dfrac{k_2(\varphi)}{k_1(\varphi)}$. Then corresponding matrix form of (5) can be given

$$\begin{bmatrix} d\vec{q}/d\varphi \\ d\vec{h}/d\varphi \\ d\vec{a}/d\varphi \end{bmatrix} = \begin{bmatrix} 0 & 1 & 0 \\ -\varepsilon & 0 & f(\varphi) \\ 0 & \varepsilon f(\varphi) & 0 \end{bmatrix} \begin{bmatrix} \vec{q} \\ \vec{h} \\ \vec{a} \end{bmatrix}. \quad (10)$$

From the first and second equations of new Frenet derivatives (10) we have

$$\vec{a} = \dfrac{1}{f(\varphi)}\left(\dfrac{d^2 \vec{q}}{d\varphi^2} + \varepsilon \vec{q}\right). \quad (11)$$

Substituting the above equation in the last equation of (10) we have the first equation of (9).

If $N$ is a spacelike ruled surface, then considering Frenet formulae (7) and following the same procedure we have the second equation of (9) immediately.

## 4. Timelike Similar Ruled Surfaces in Minkowski 3-space $E_1^3$

In this section we introduce the definition and characterizations of timelike similar ruled surfaces with variable transformation in $E_1^3$. First, we give the following definition.

**Definition 4.1.** Let $N_\alpha$ and $N_\beta$ be two regular timelike ruled surfaces of the same type in $E_1^3$ given by the parametrizations

$$\begin{cases} \vec{r}_\alpha(s_\alpha, v_\alpha) = \vec{\alpha}(s_\alpha) + v_\alpha \vec{q}_\alpha(s_\alpha), \\ \vec{r}_\beta(s_\beta, v_\beta) = \vec{\beta}(s_\beta) + v_\beta \vec{q}_\beta(s_\beta), \end{cases} \quad (12)$$

respectively, where $\vec{\alpha}(s_\alpha)$ and $\vec{\beta}(s_\beta)$ are striction curves of $N_\alpha$ and $N_\beta$ and $s_\alpha$, $s_\beta$ are arc length parameters of $\vec{\alpha}(s_\alpha)$ and $\vec{\beta}(s_\beta)$, respectively. Let the Frenet frames and invariants of

$N_\alpha$ and $N_\beta$ be $\{\vec{q}_\alpha, \vec{h}_\alpha, \vec{a}_\alpha\}$, $k_1^\alpha, k_2^\alpha$ and $\{\vec{q}_\beta, \vec{h}_\beta, \vec{a}_\beta\}$, $k_1^\beta, k_2^\beta$, respectively. $N_\alpha$ and $N_\beta$ are called timelike similar ruled surfaces with variable transformation $\lambda_\beta^\alpha$ if there exists a variable transformation

$$s_\alpha = \int \lambda_\beta^\alpha(s_\beta) ds_\beta, \tag{13}$$

of the arc lengths such that the rulings are the same for two ruled surfaces i.e.,

$$\vec{q}_\alpha(s_\alpha) = \vec{q}_\beta(s_\beta), \tag{14}$$

for all corresponding values of parameters under the transformation $\lambda_\beta^\alpha$. All timelike ruled surfaces satisfying equation (14) are called a family of timelike similar ruled surfaces with variable transformation.

Then we can give the following theorems characterizing timelike similar ruled surfaces. Whenever we talk about $N_\alpha$ and $N_\beta$ we mean that the surfaces are regular and have the parametrizations as given in (12).

**Theorem 4.1.** *Let $N_\alpha$ and $N_\beta$ be two timelike ruled surfaces in $E_1^3$. Then $N_\alpha$ and $N_\beta$ are timelike similar ruled surfaces with variable transformation if and only if the central normal vectors of the surfaces are the same, i.e.,*

$$\vec{h}_\alpha(s_\alpha) = \vec{h}_\beta(s_\beta), \tag{15}$$

*under the particular variable transformation*

$$\lambda_\beta^\alpha = \frac{ds_\alpha}{ds_\beta} = \frac{k_1^\beta}{k_1^\alpha}, \tag{16}$$

*of the arc lengths.*

**Proof:** Let $N_\alpha$ and $N_\beta$ be two timelike similar ruled surfaces in $E_1^3$ with variable transformation. Then differentiating (14) with respect to $s_\beta$ it follows

$$k_1^\alpha \lambda_\beta^\alpha \vec{h}_\alpha = k_1^\beta \vec{h}_\beta. \tag{17}$$

From (17) we obtain (15) and (16) immediately.

Conversely, let $N_\alpha$ and $N_\beta$ be two regular timelike ruled surfaces in $E_1^3$ satisfying (15) and (16). By multiplying (15) with $k_1^\beta$ and differentiating the results equality with respect to $s_\beta$ we have

$$\int k_1^\beta(s_\beta)\vec{h}_\beta(s_\beta)ds_\beta = \int k_1^\beta(s_\beta)\vec{h}_\beta(s_\beta)\frac{ds_\beta}{ds_\alpha}ds_\alpha. \tag{18}$$

From (15) and (16) we obtain

$$\vec{q}_\beta(s_\beta) = \int k_1^\beta(s_\beta)\vec{h}_\beta(s_\beta)ds_\beta = \int k_1^\alpha(s_\alpha)\vec{h}_\alpha(s_\alpha)ds_\alpha = \vec{q}_\alpha(s_\alpha), \tag{19}$$

which means that $N_\alpha$ and $N_\beta$ are timelike similar ruled surfaces with variable transformation.

**Theorem 4.2.** *Let $N_\alpha$ and $N_\beta$ be two timelike ruled surfaces in $E_1^3$. Then $N_\alpha$ and $N_\beta$ are timelike similar ruled surfaces with variable transformation if and only if the asymptotic normal vectors of the surfaces are satisfied the following equality*

$$\vec{a}_\alpha(s_\alpha) = \varepsilon_\alpha \varepsilon_\beta \vec{a}_\beta(s_\beta), \tag{20}$$

*under the particular variable transformation*

$$\lambda_\beta^\alpha = \frac{ds_\alpha}{ds_\beta} = \frac{k_2^\beta}{k_2^\alpha}, \tag{21}$$

*of the arc lengths, where* $\varepsilon_\alpha = \langle \vec{q}_\alpha, \vec{q}_\alpha \rangle = \pm 1$, $\varepsilon_\beta = \langle \vec{q}_\beta, \vec{q}_\beta \rangle = \pm 1$.

**Proof:** Let $N_\alpha$ and $N_\beta$ be two timelike similar ruled surfaces in $E_1^3$ with variable transformation. Then from Definition 4.1 and Theorem 4.1 there exists a variable transformation of the arc lengths such that the rulings and central normal vectors are the same. Then from (14) and (15) we have

$$\vec{a}_\alpha(s_\alpha) = \varepsilon_\alpha \left( \vec{q}_\alpha(s_\alpha) \times \vec{h}_\alpha(s_\alpha) \right) = \varepsilon_\alpha \left( \vec{q}_\beta(s_\beta) \times \vec{h}_\beta(s_\beta) \right) = \varepsilon_\alpha \varepsilon_\beta \vec{a}_\beta(s_\beta). \tag{22}$$

Conversely, let $N_\alpha$ and $N_\beta$ be two timelike ruled surfaces in $E_1^3$ satisfying (20) and (21). By differentiating (20) with respect to $s_\beta$ we obtain

$$\varepsilon_\alpha k_2^\alpha(s_\alpha) \vec{h}_\alpha(s_\alpha) \frac{ds_\alpha}{ds_\beta} = \varepsilon_\alpha \varepsilon_\beta \left( \varepsilon_\beta k_2^\beta(s_\beta) \vec{h}_\beta(s_\beta) \right), \tag{23}$$

which gives us

$$\lambda_\beta^\alpha = \frac{k_2^\beta}{k_2^\alpha}, \quad \vec{h}_\alpha(s_\alpha) = \vec{h}_\beta(s_\beta). \tag{24}$$

Then from (20) and (24) we have

$$\begin{aligned}\vec{q}_\alpha(s_\alpha) &= -\varepsilon_\alpha \vec{h}_\alpha(s_\alpha) \times \vec{a}_\alpha(s_\alpha) = -\varepsilon_\alpha \left( \varepsilon_\alpha \varepsilon_\beta \vec{h}_\beta(s_\beta) \times \vec{a}_\beta(s_\beta) \right) = -\varepsilon_\beta \vec{h}_\beta(s_\beta) \times \vec{a}_\beta(s_\beta) \\ &= \vec{q}_\beta(s_\beta) \end{aligned} \tag{25}$$

which completes the proof.

***Theorem 4.3.*** *Let $N_\alpha$ and $N_\beta$ be two timelike ruled surfaces in $E_1^3$. Then $N_\alpha$ and $N_\beta$ are timelike similar ruled surfaces with variable transformation if and only if the ratio of curvatures are the same i.e.,*

$$\frac{k_2^\beta(s_\beta)}{k_1^\beta(s_\beta)} = \frac{k_2^\alpha(s_\alpha)}{k_1^\alpha(s_\alpha)}, \tag{26}$$

*under the particular variable transformation keeping equal total curvatures, i.e.,*

$$\varphi_\beta(s_\beta) = \int k_1^\beta(s_\beta) ds_\beta = \int k_1^\alpha(s_\alpha) ds_\alpha = \varphi_\alpha(s_\alpha) \tag{27}$$

*of the arc lengths.*

**Proof:** Let $N_\alpha$ and $N_\beta$ be two timelike similar ruled surfaces in $E_1^3$ with variable transformation. Then from (21) and (24) we have (26) under the variable transformation (27), and this transformation is also leads from (21) by integration.

Conversely, let $N_\alpha$ and $N_\beta$ be two timelike ruled surfaces in $E_1^3$ satisfying (26) and (27). From Theorem 3.2, the rulings $\vec{q}_\alpha$ and $\vec{q}_\beta$ of the surfaces $N_\alpha$ and $N_\beta$ satisfy the following vector differential equations of third order

$$\frac{d}{d\varphi_\alpha}\left( \frac{1}{f_\alpha(\varphi_\alpha)} \frac{d^2\vec{q}_\alpha}{d\varphi_\alpha^2} \right) + \varepsilon_\alpha \left( \frac{1 - f_\alpha^2(\varphi_\alpha)}{f_\alpha(\varphi_\alpha)} \right) \frac{d\vec{q}_\alpha}{d\varphi_\alpha} - \varepsilon_\alpha \left( \frac{1}{f_\alpha^2(\varphi_\alpha)} \frac{df_\alpha(\varphi_\alpha)}{d\varphi_\alpha} \right) \vec{q}_\alpha = 0, \tag{28}$$

$$\frac{d}{d\varphi_\beta}\left(\frac{1}{f_\beta(\varphi_\beta)}\frac{d^2\vec{q}_\beta}{d\varphi_\beta^2}\right)+\varepsilon_\beta\left(\frac{1-f_\beta^2(\varphi_\beta)}{f_\beta(\varphi_\beta)}\right)\frac{d\vec{q}_\beta}{d\varphi_\beta}-\varepsilon_\beta\left(\frac{1}{f_\beta^2(\varphi_\beta)}\frac{df_\beta(\varphi_\beta)}{d\varphi_\beta}\right)\vec{q}_\beta=0, \quad (29)$$

respectively, where

$$f_\alpha(\varphi_\alpha)=\frac{k_2^\alpha(\varphi_\alpha)}{k_1^\alpha(\varphi_\alpha)},\quad f_\beta(\varphi_\beta)=\frac{k_2^\beta(\varphi_\beta)}{k_1^\beta(\varphi_\beta)},\quad \varphi_\alpha(s_\alpha)=\int k_1^\alpha(s_\alpha)ds_\alpha,\quad \varphi_\beta(s_\beta)=\int k_1^\beta(s_\beta)ds_\beta.$$

From (26) we have $f_\alpha(\varphi_\alpha)=f_\beta(\varphi_\beta)$ under the variable transformation $\varphi_\alpha=\varphi_\beta$. Thus under the equation (26) and transformation (27), the equations (28) and (29) are the same, i.e., they have the same solutions. It means that the rulings $\vec{q}_\alpha$ and $\vec{q}_\beta$ are the same. Then $N_\alpha$ and $N_\beta$ are two timelike similar ruled surfaces in $E_1^3$ with variable transformation.

**Theorem 4.4.** *Let timelike ruled surfaces $N_\alpha$ and $N_\beta$ be developable surfaces and let the striction lines have the same Lorentzian characters with the rulings. Then $N_\alpha$ and $N_\beta$ are timelike similar ruled surfaces with variable transformation if and only if the striction curves of the surfaces are similar curves with variable transformation.*

**Proof:** Let developable timelike ruled surfaces $N_\alpha$ and $N_\beta$ be two timelike similar ruled surfaces in $E_1^3$ with variable transformation. Since the surfaces are developable, from Theorem 3.1 we have

$$\frac{d\vec{\alpha}}{ds_\alpha}=\vec{T}_\alpha(s_\alpha)=\vec{q}_\alpha(s_\alpha),\quad \frac{d\vec{\beta}}{ds_\beta}=\vec{T}_\beta(s_\beta)=\vec{q}_\beta(s_\beta). \quad (30)$$

where $\vec{T}_\alpha(s_\alpha)$ and $\vec{T}_\beta(s_\beta)$ are unit tangents of the striction curves $\vec{\alpha}(s_\alpha)$ and $\vec{\beta}(s_\beta)$, respectively. From (14) and (30) we have

$$\frac{d\vec{\alpha}}{ds_\alpha}=\vec{q}_\alpha(s_\alpha)=\vec{q}_\beta(s_\beta)=\frac{d\vec{\beta}}{ds_\beta}, \quad (31)$$

which shows that striction curves $\vec{\alpha}(s_\alpha)$ and $\vec{\beta}(s_\beta)$ are similar curves in $E_1^3$.

Conversely, if the striction curves $\vec{\alpha}(s_\alpha)$ and $\vec{\beta}(s_\beta)$ are similar curves, then there exists a variable transformation between arc lengths such that

$$\frac{d\vec{\alpha}}{ds_\alpha}=\vec{T}_\alpha(s_\alpha)=\vec{T}_\beta(s_\beta)=\frac{d\vec{\beta}}{ds_\beta}. \quad (32)$$

Since the ruled surfaces are developable, from Theorem 3.1 we have $\vec{T}_\alpha(s_\alpha)=\vec{q}_\alpha(s_\alpha)$ and $\vec{T}_\beta(s_\beta)=\vec{q}_\beta(s_\beta)$. Then from (32) we have that $\vec{q}_\alpha(s_\alpha)=\vec{q}_\beta(s_\beta)$, i.e., $N_\alpha$ and $N_\beta$ are timelike similar ruled surfaces with variable transformation.

Let now consider some special cases. From (16) and (24) we have
$$k_1^\beta=\lambda_\beta^\alpha k_1^\alpha,\quad k_2^\beta=\lambda_\beta^\alpha k_2^\alpha, \quad (33)$$
respectively. From (33) it is clear that if $N_\alpha$ is a timelike cylindrical surface i.e., $k_1^\alpha=0$, then under the variable transformation the curvature does not change. So we have the following corollaries.

***Corollary 4.1.*** *The family of timelike cylindrical surfaces forms a family of timelike similar ruled surfaces with variable transformation.*

If $N_\alpha$ is a timelike conoid surface i.e., $k_2^\alpha = 0$, then under the variable transformation the curvature does not change. So we have the following corollary.

***Corollary 4.2.*** *The family of timelike conoid surfaces forms a family of timelike similar ruled surfaces with variable transformation.*

**5. Spacelike Similar Ruled Surfaces in Minkowski 3-space $E_1^3$**

In this section we introduce the definition and characterizations of spacelike similar ruled surfaces with variable transformation in $E_1^3$. First, we give the following definition.

**Definition 5.1.** Let $N_\alpha$ and $N_\beta$ be two regular spacelike ruled surfaces in $E_1^3$ given by the parametrizations
$$\begin{cases} \vec{r}_\alpha(s_\alpha, v_\alpha) = \vec{\alpha}(s_\alpha) + v_\alpha \vec{q}_\alpha(s_\alpha), \\ \vec{r}_\beta(s_\beta, v_\beta) = \vec{\beta}(s_\beta) + v_\beta \vec{q}_\beta(s_\beta), \end{cases} \quad (34)$$
respectively, where $\vec{\alpha}(s_\alpha)$ and $\vec{\beta}(s_\beta)$ are striction curves of $N_\alpha$ and $N_\beta$ and $s_\alpha$, $s_\beta$ are arc length parameters of $\vec{\alpha}(s_\alpha)$ and $\vec{\beta}(s_\beta)$, respectively. Let the Frenet frames and invariants of $N_\alpha$ and $N_\beta$ be $\{\vec{q}_\alpha, \vec{h}_\alpha, \vec{a}_\alpha\}$, $k_1^\alpha, k_2^\alpha$ and $\{\vec{q}_\beta, \vec{h}_\beta, \vec{a}_\beta\}$, $k_1^\beta, k_2^\beta$, respectively. $N_\alpha$ and $N_\beta$ are called spacelike similar ruled surfaces with variable transformation $\lambda_\beta^\alpha$ if there exists a variable transformation
$$s_\alpha = \int \lambda_\beta^\alpha(s_\beta) ds_\beta, \quad (35)$$
of the arc lengths such that the rulings are the same for two ruled surfaces i.e.,
$$\vec{q}_\alpha(s_\alpha) = \vec{q}_\beta(s_\beta), \quad (36)$$
for all corresponding values of parameters under the transformation $\lambda_\beta^\alpha$. All spacelike ruled surfaces satisfying equation (36) are called a family of spacelike similar ruled surfaces with variable transformation.

Then we can give the following theorems characterizing spacelike similar ruled surfaces. Whenever we talk about $N_\alpha$ and $N_\beta$ we mean that the surfaces are regular and have the parametrizations as given in (34).

***Theorem 5.1.*** *Let $N_\alpha$ and $N_\beta$ be two spacelike ruled surfaces in $E_1^3$. Then $N_\alpha$ and $N_\beta$ are spacelike similar ruled surfaces with variable transformation if and only if the central normal vectors of the surfaces are the same, i.e.,*
$$\vec{h}_\alpha(s_\alpha) = \vec{h}_\beta(s_\beta), \quad (37)$$
*under the particular variable transformation*
$$\lambda_\beta^\alpha = \frac{ds_\alpha}{ds_\beta} = \frac{k_1^\beta}{k_1^\alpha}, \quad (38)$$
*of the arc lengths.*

**Proof:** Let $N_\alpha$ and $N_\beta$ be two spacelike similar ruled surfaces in $E_1^3$ with variable transformation. Then differentiating (36) with respect to $s_\beta$ it follows

$$k_1^\alpha \lambda_\beta^\alpha \vec{h}_\alpha = k_1^\beta \vec{h}_\beta. \tag{39}$$

From (39) we obtain (37) and (38) immediately.

Conversely, let $N_\alpha$ and $N_\beta$ be two spacelike ruled surfaces in $E_1^3$ satisfying (37) and (38). By multiplying (37) with $k_1^\beta$ and differentiating the results equality with respect to $s_\beta$ we have

$$\int k_1^\beta(s_\beta)\vec{h}_\beta(s_\beta)ds_\beta = \int k_1^\beta(s_\beta)\vec{h}_\beta(s_\beta)\frac{ds_\beta}{ds_\alpha}ds_\alpha. \tag{40}$$

From (37) and (38) we obtain

$$\vec{q}_\beta(s_\beta) = \int k_1^\beta(s_\beta)\vec{h}_\beta(s_\beta)ds_\beta = \int k_1^\alpha(s_\alpha)\vec{h}_\alpha(s_\alpha)ds_\alpha = \vec{q}_\alpha(s_\alpha), \tag{41}$$

which means that $N_\alpha$ and $N_\beta$ are spacelike similar ruled surfaces with variable transformation.

**Theorem 5.2.** *Let $N_\alpha$ and $N_\beta$ be two spacelike ruled surfaces in $E_1^3$. Then $N_\alpha$ and $N_\beta$ are spacelike similar ruled surfaces with variable transformation if and only if the asymptotic normal vectors of the surfaces are the same i.e.,*

$$\vec{a}_\alpha(s_\alpha) = \vec{a}_\beta(s_\beta), \tag{42}$$

*under the particular variable transformation*

$$\lambda_\beta^\alpha = \frac{ds_\alpha}{ds_\beta} = \frac{k_2^\beta}{k_2^\alpha}, \tag{43}$$

*of the arc lengths.*

**Proof:** Let $N_\alpha$ and $N_\beta$ be two spacelike similar ruled surfaces in $E_1^3$ with variable transformation. Then from Definition 5.1 and Theorem 5.1 there exists a variable transformation of the arc lengths such that the rulings and central normal vectors are the same. Then from (36) and (37) we have

$$\vec{a}_\alpha(s_\alpha) = -\vec{q}_\alpha(s_\alpha) \times \vec{h}_\alpha(s_\alpha) = -\vec{q}_\beta(s_\beta) \times \vec{h}_\beta(s_\beta) = \vec{a}_\beta(s_\beta). \tag{44}$$

Conversely, let $N_\alpha$ and $N_\beta$ be two spacelike ruled surfaces in $E_1^3$ satisfying (42) and (43). By differentiating (42) with respect to $s_\beta$ we obtain

$$k_2^\alpha(s_\alpha)\vec{h}_\alpha(s_\alpha)\frac{ds_\alpha}{ds_\beta} = k_2^\beta(s_\beta)\vec{h}_\beta(s_\beta), \tag{45}$$

which gives us

$$\lambda_\beta^\alpha = \frac{k_2^\beta}{k_2^\alpha}, \quad \vec{h}_\alpha(s_\alpha) = \vec{h}_\beta(s_\beta). \tag{46}$$

Then from (42) and (46) we have

$$\vec{q}_\alpha(s_\alpha) = -\vec{h}_\alpha(s_\alpha) \times \vec{a}_\alpha(s_\alpha) = -\vec{h}_\beta(s_\beta) \times \vec{a}_\beta(s_\beta) = \vec{q}_\beta(s_\beta), \tag{47}$$

which completes the proof.

**Theorem 5.3.** *Let $N_\alpha$ and $N_\beta$ be two spacelike ruled surfaces in $E_1^3$. Then $N_\alpha$ and $N_\beta$ are spacelike similar ruled surfaces with variable transformation if and only if the ratio of curvatures are the same i.e.,*

$$\frac{k_2^\beta(s_\beta)}{k_1^\beta(s_\beta)} = \frac{k_2^\alpha(s_\alpha)}{k_1^\alpha(s_\alpha)}, \tag{48}$$

*under the particular variable transformation keeping equal total curvatures, i.e.,*

$$\varphi_\beta(s_\beta) = \int k_1^\beta(s_\beta)ds_\beta = \int k_1^\alpha(s_\alpha)ds_\alpha = \varphi_\alpha(s_\alpha) \tag{49}$$

*of the arc lengths.*

**Proof:** Let $N_\alpha$ and $N_\beta$ be two spacelike similar ruled surfaces in $E_1^3$ with variable transformation. Then from (43) and (46) we have (48) under the variable transformation (49), and this transformation is also leads from (43) by integration.

Conversely, let $N_\alpha$ and $N_\beta$ be two regular spacelike ruled surfaces in $E_1^3$ satisfying (48) and (49). From Theorem 3.2, the rulings $\vec{q}_\alpha$ and $\vec{q}_\beta$ of the surfaces $N_\alpha$ and $N_\beta$ satisfy the following vector differential equations of third order

$$\frac{d}{d\varphi_\alpha}\left(\frac{1}{f_\alpha(\varphi_\alpha)}\frac{d^2\vec{q}_\alpha}{d\varphi_\alpha^2}\right) - \left(\frac{1+f_\alpha^2(\varphi_\alpha)}{f_\alpha(\varphi_\alpha)}\right)\frac{d\vec{q}_\alpha}{d\varphi_\alpha} + \left(\frac{1}{f_\alpha^2(\varphi_\alpha)}\frac{df_\alpha(\varphi_\alpha)}{d\varphi_\alpha}\right)\vec{q}_\alpha = 0, \tag{50}$$

$$\frac{d}{d\varphi_\beta}\left(\frac{1}{f_\beta(\varphi_\beta)}\frac{d^2\vec{q}_\beta}{d\varphi_\beta^2}\right) - \left(\frac{1+f_\beta^2(\varphi_\beta)}{f_\beta(\varphi_\beta)}\right)\frac{d\vec{q}_\beta}{d\varphi_\beta} + \left(\frac{1}{f_\beta^2(\varphi_\beta)}\frac{df_\beta(\varphi_\beta)}{d\varphi_\beta}\right)\vec{q}_\beta = 0, \tag{51}$$

where

$$f_\alpha(\varphi_\alpha) = \frac{k_2^\alpha(\varphi_\alpha)}{k_1^\alpha(\varphi_\alpha)}, \quad f_\beta(\varphi_\beta) = \frac{k_2^\beta(\varphi_\beta)}{k_1^\beta(\varphi_\beta)}, \quad \varphi_\alpha(s_\alpha) = \int k_1^\alpha(s_\alpha)ds_\alpha, \quad \varphi_\beta(s_\beta) = \int k_1^\beta(s_\beta)ds_\beta.$$

From (26) we have $f_\alpha(\varphi_\alpha) = f_\beta(\varphi_\beta)$ under the variable transformation $\varphi_\alpha = \varphi_\beta$. Under the equation (48) and transformation (49), the equations (50) and (51) are the same, i.e., they have the same solutions. It means that the rulings $\vec{q}_\alpha$ and $\vec{q}_\beta$ are the same. Then $N_\alpha$ and $N_\beta$ are two spacelike similar ruled surfaces in $E_1^3$ with variable transformation.

**Theorem 5.4.** *Let spacelike ruled surfaces $N_\alpha$ and $N_\beta$ be developable surfaces. Then $N_\alpha$ and $N_\beta$ are spacelike similar ruled surfaces with variable transformation if and only if the striction curves of the surfaces are similar curves with variable transformation.*

**Proof:** Let developable spacelike ruled surfaces $N_\alpha$ and $N_\beta$ be two spacelike similar ruled surfaces in $E_1^3$ with variable transformation. Since the surfaces are developable, from Theorem 3.1 we have

$$\frac{d\vec{\alpha}}{ds_\alpha} = \vec{T}_\alpha(s_\alpha) = \vec{q}_\alpha(s_\alpha), \quad \frac{d\vec{\beta}}{ds_\beta} = \vec{T}_\beta(s_\beta) = \vec{q}_\beta(s_\beta), \tag{52}$$

where $\vec{T}_\alpha(s_\alpha)$ and $\vec{T}_\beta(s_\beta)$ are unit tangents of the striction curves $\vec{\alpha}(s_\alpha)$ and $\vec{\beta}(s_\beta)$, respectively. From (36) and (52) we have

$$\frac{d\vec{\alpha}}{ds_\alpha} = \vec{q}_\alpha(s_\alpha) = \vec{q}_\beta(s_\beta) = \frac{d\vec{\beta}}{ds_\beta} \tag{53}$$

which shows that striction curves $\vec{\alpha}(s_\alpha)$ and $\vec{\beta}(s_\beta)$ are similar curves in $E_1^3$.

Conversely, if the striction curves $\vec{\alpha}(s_\alpha)$ and $\vec{\beta}(s_\beta)$ are similar curves, then there exists a variable transformation between arc lengths such that

$$\frac{d\vec{\alpha}}{ds_\alpha} = \vec{T}_\alpha(s_\alpha) = \vec{T}_\beta(s_\beta) = \frac{d\vec{\beta}}{ds_\beta}. \tag{54}$$

Since the ruled surfaces are developable, from Theorem 3.1 we have $\vec{T}_\alpha(s_\alpha) = \vec{q}_\alpha(s_\alpha)$ and $\vec{T}_\beta(s_\beta) = \vec{q}_\beta(s_\beta)$. From (54) we have that $\vec{q}_\alpha(s_\alpha) = \vec{q}_\beta(s_\beta)$, i.e., $N_\alpha$ and $N_\beta$ are spacelike similar ruled surfaces with variable transformation.

Let now consider some special cases. From (38) and (46) we have
$$k_1^\beta = \lambda_\beta^\alpha k_1^\alpha, \quad k_2^\beta = \lambda_\beta^\alpha k_2^\alpha, \tag{55}$$
respectively. From (55) it is clear that if $N_\alpha$ is a cylindrical surface i.e., $k_1^\alpha = 0$, then under the variable transformation the curvature does not change. So we have the following corollaries.

***Corollary 5.1.*** *The family of spacelike cylindrical surfaces forms a family of spacelike similar ruled surfaces with variable transformation.*

If $N_\alpha$ is a spacelike conoid surface i.e., $k_2^\alpha = 0$, then under the variable transformation the curvature does not change. So we have the following corollary.

***Corollary 5.2.*** *The family of spacelike conoid surfaces forms a family of spacelike similar ruled surfaces with variable transformation.*

## 6. Conclusions

In Minkowski 3-space, some special families of timelike and spacelike ruled surfaces are defined and called similar ruled surfaces. Some properties of these special surfaces are obtained and it is showed that developable ruled surfaces form a family of similar ruled surfaces in $E_1^3$ if and only if the striction curves of the surfaces are similar curves with variable transformation in $E_1^3$. Of course, in Minkowski 3-space another type of the ruled surfaces is ruled surface with lightlike ruling. One can consider the similar ruled surfaces with lightlike ruling by considering this present paper and can obtain similar characterizations for these surfaces.